# A Dedekind Psi Function Inequality
N. A. Carella, December, 2011.

***Abstract***: This note shows that the arithmetic function $\psi(N)/N = \prod_{p|N}(1+1/p)$, called the Dedekind psi function, achieves its extreme values on the subset of primorial integers $N = 2^{v_1} \cdot 3^{v_2} \cdots p_k^{v_k}$, where $p_i$ is the $k$th prime, and $v_i \geq 1$. In particular, the inequality $\psi(N)/N > 6\pi^{-2}e^{\gamma} \log\log N$ holds for all large squarefree primorial integers $N = 2 \cdot 3 \cdot 5 \cdots p_k$ unconditionally.

## 1 Introduction

The psi function $\psi(N) = N\prod_{p|N}(1+1/p)$ and its normalized counterpart $\psi(N)/N = \prod_{p|N}(1+1/p)$ arise in various mathematic, and physic problems. Moreover, this function is entangled with other arithmetic functions. The values of the normalized psi function coincide with the squarefree kernel

$$\sum_{d|N} \frac{\mu^2(d)}{d} \qquad (1)$$

of the sum of divisor function $\sigma(N) = \sum_{d|N} d$. In particular, $\psi(N)/N = \prod_{p|N}(1+1/p) = \sigma(N)/N$ on the subset of square-free integers. This note proposes a new lower estimate of the Dedekind function.

***Theorem 1.*** Let $N \in \mathbb{N}$ be an integer, then $\psi(N)/N > 6\pi^{-2}e^{\gamma} \log\log N$ holds unconditionally for all sufficiently large primorial integer $N = 2\cdot 3\cdot 5\cdots p_k$, where $p_k$ is the $k$th prime.

An intuitive and clear-cut relationship between the Riemann hypothesis and the Dedekind psi function is established in Theorem 4 by means of the prime number theorem. Recall that the Riemann hypothesis claims that the nontrivial zeros of the zeta function ζ(s) are located on the critical line { $\Re\mathrm{e}(s) = 1/2$ }, see [ES]. A survey of various strikingly different reformulations of the Riemann hypothesis appears in [AM]. And a derivation of the reformulation of the Riemann hypothesis in terms of the psi function was recently resolved in [SP].

This estimate is quite similar to the lower estimate of the totient function, namely,

$$N_k / \varphi(N_k) > e^{\gamma} \log\log N_k \tag{2}$$

for all large squarefree primorial integers $N_k = 2 \cdot 3 \cdot 5 \cdots p_k$. This inequality, called the Nicolas inequality, is also a reformulation of the Riemann hypothesis, see [NS].

The next section discusses background information, and supplies the proof of Theorem 1 in Subsection 2.2. The derivation of this result unfolds from a result on the oscillations theorem of finite prime product.

## 2 Elementary Materials And The Proof
The basic concepts and results employed throughout this work are stated in this Section.

**2.1 Sums and Products Over the Primes.** The most basic finite sum over the prime numbers is the prime harmonic sum $\sum_{n \leq x} p^{-1}$. The refined estimate of this finite sum, stated below, is a synthesis of various results due to various authors. The earliest version $\sum_{n \leq x} p^{-1} = \log\log x + B_1 + O(1/\log x)$ is due to Mertens.

***Theorem 2.*** Let $x \geq 2$ be a sufficiently large number. Then

$$\sum_{p \leq x} \frac{1}{p} = \begin{cases} \log\log x + B_1 + O(e^{-c(\log\log x)^{1/2}}), & \text{unconditionally,} \\ \log\log x + B_1 + O((\log x)^{-1/2}), & \text{conditional on the Riemann hypothesis,} \\ \log\log x + B_1 + \Omega_{\pm}((\log x)^{-1/2} \log\log\log x / \log\log x), & \text{unconditional oscillations,} \end{cases} \tag{3}$$

where $B_1 = .2614972128\ldots$.

Proof: Use the integral representation of the finite sum

$$\sum_{p \leq x} \frac{1}{p} = \int_c^x \frac{d\pi(t)}{t}, \tag{4}$$

where $c > 1$ is a small constant. Here, the prime counting function $\pi(x) = \#\{ p \leq x : p \text{ is prime} \}$ has the form

$$\pi(x) = \begin{cases} li(x) + O(xe^{-c(\log x)^{1/2}}), & \text{unconditionally,} \\ li(x) + O(x^{1/2} \log x), & \text{conditional on the Riemann hypothesis,} \\ li(x) + \Omega_{\pm}(x^{1/2} \log\log\log x / \log x), & \text{unconditional oscillations.} \end{cases} \tag{5}$$

The unconditional part of the prime counting formula arises from the delaVallee Poussin form of the prime number theorem $\pi(x) = li(x) + O(xe^{-c(\log x)^{1/2}})$, see [MV, p. 179], the conditional part arises from the Riemann form of the prime number theorem $\pi(x) = li(x) + O(x^{1/2} \log x)$, and the unconditional oscillations part arises from the Littlewood form of the prime number theorem $\pi(x) = li(x) + \Omega_{\pm}(x^{1/2} \log\log\log x / \log x)$, consult [IV, p.

51], [MV, p. 479] et cetera. Now proceed to replace the logarithm integral $li(x) = \int_0^x (t \log t)^{-1} dt$, and the appropriate prime counting measure $d\pi(t)$, and simplify the integral. ∎

The proof of the unconditional part of this result is widely available in the literature, see [HW], [MV], [TN], et cetera. As an application of the last result, there are the following interesting products:

The omega notation $f(x) = g(x) + \Omega_\pm(h(x))$ means that both $f(x) > g(x) + c_0 h(x)$ and $f(x) < g(x) - c_0 h(x)$ occur infinitely often as $x \to \infty$, where $c_0 > 0$ is a constant, see [MV, p. 5], [WK].

**Theorem 3.** Let $x \in \mathbb{R}$ be a large real number, then

$$\prod_{p \leq x}\left(1 - \frac{1}{p}\right)^{-1} = \begin{cases} e^\gamma \log x + O(e^{-c(\log x)^{1/2}} \log x), & \text{unconditionally}, \\ e^\gamma \log x + O(x^{-1/2} \log x), & \text{conditional on the Riemann hypothesis}, \\ e^\gamma \log x + \Omega_\pm(x^{-1/2} \log\log\log x / \log x), & \text{unconditional oscillations}, \end{cases} \quad (6)$$

Proof: Consider the logarithm of the product

$$\log \prod_{p \leq x}(1 - 1/p)^{-1} = \sum_{p \leq x} \frac{1}{p} + \sum_{p \leq x}\sum_{n \geq 2} \frac{1}{np^n} = \sum_{p \leq x} \frac{1}{p} + \gamma - B_1 + O(1/x), \quad (7)$$

where the Euler constant is defined by $\gamma = \lim_{x \to \infty} \sum_{n \leq x}(n^{-1} - \log n) = 0.577215665...$, and the Mertens constant is defined by $B_1 = \gamma + \sum_{p \geq 2}(\log(1 - 1/p) + 1/p) = .2614972128...$, see [HW, p. 466]. The last equality in (7) stems from the power series expansion $B_1 = \gamma - \sum_{p \geq 2}\sum_{n \geq 2}(np^n)^{-1}$, which yields

$$\sum_{p \leq x}\sum_{n \geq 2}\frac{1}{np^n} = \gamma - B_1 - \sum_{p > x}\sum_{n \geq 2}\frac{1}{np^n} = \gamma - B_1 + O(1/x), \quad (8)$$

The remaining steps follows from Theorem 2, and reversing the logarithm. ∎

The third part above simplifies the proof given in [DP] of the following result:

The quantity

$$x^{1/2}\left(\prod_{p \leq x}(1 - 1/p)^{-1} - e^\gamma \log x\right) \quad (9)$$

attains arbitrary large positive and negative values as $x \to \infty$.

**Theorem 4.** Let $x \geq x_0$ be a real number, then

$$\prod_{p \leq x}\left(1+\frac{1}{p}\right) = \begin{cases} 6\pi^{-2}e^{\gamma}\log x + O(e^{-c(\log x)^{1/2}}\log x), & \text{unconditionally,} \\ 6\pi^{-2}e^{\gamma}\log x + O(x^{-1/2}\log x), & \text{conditional on the Riemann hypothesis,} \\ 6\pi^{-2}e^{\gamma}\log x + \Omega_{\pm}(x^{-1/2}\log\log x / \log x), & \text{unconditional oscillations,} \end{cases} \quad (10)$$

Proof: For a large real number $x \in \mathbb{R}$, rewrite the product as

$$\prod_{p \leq x}(1+1/p)(1-1/p)(1-1/p)^{-1} = \prod_{p \leq x}(1-1/p^2)\prod_{p \leq x}(1-1/p)^{-1}. \quad (11)$$

Replacing $\prod_{p \geq 2}(1-1/p^2) = \sum_{n \geq 1}\mu(n)n^{-2} = 6\pi^{-2}$ in the first product on the right side, yields

$$\begin{aligned}\prod_{p \leq x}(1-1/p^2)\prod_{p \leq x}(1-1/p)^{-1} &= \left(\frac{6}{\pi^2} - \sum_{n \geq z}\mu(n)n^{-2}\right)\prod_{p \leq x}(1-1/p)^{-1} \\ &= \frac{6}{\pi^2}\prod_{p \leq x}(1-1/p)^{-1} + O(\frac{\log x}{x^c}),\end{aligned} \quad (12)$$

where $z = O(x^c)$, $c \geq 1$ constant. Lastly, applying Theorem 3, to the last product above, yields the claim. ∎

This result immediately gives an improved upper bound of the sum of divisors function; currently, the best upper bound of the sum of divisors function is

$$\sigma(N)/N \leq e^{\gamma}\log\log N + O(1/\log\log N) \quad (13)$$

for any integer $N \geq 1$, see [RN].

**Corollary 5.** The sum of divisors function satisfies the followings inequalities unconditionally:

(i) $\displaystyle\sum_{d \mid N}\frac{1}{d} = \begin{cases} e^{\gamma}\log\log N + O(e^{-c(\log\log N)^{1/2}}), & \text{unconditionally,} \\ e^{\gamma}\log\log N + O((\log N)^{-1/2}\log\log N), & \text{conditional on the RH,} \\ e^{\gamma}\log\log N + \Omega_{\pm}((\log N)^{-1/2}\log\log\log N / \log\log N), & \text{unconditional oscillations,} \end{cases}$

the sum of divisors function $\sigma(N)/N$, (14)

(ii) $\displaystyle\sum_{d \mid N, \mu(d) \neq 0}\frac{1}{d} = \begin{cases} 6\pi^{-1}e^{\gamma}\log\log N + O(e^{-c(\log\log N)^{1/2}}), & \text{unconditionally,} \\ 6\pi^{-1}e^{\gamma}\log\log N + O((\log N)^{-1/2}\log\log N), & \text{conditional on the RH,} \\ 6\pi^{-1}e^{\gamma}\log\log N + \Omega_{\pm}((\log N)^{-1/2}\log\log\log N / \log\log N), & \text{unconditional oscillations,} \end{cases}$

the squarefree kernel of $\sigma(N)/N$,

(ii) $\displaystyle\sum_{d|N,\,\mu(d)\neq 0}\frac{1}{d} = \begin{cases} (1-6\pi^{-1})e^{\gamma}\log\log N + O(e^{-c(\log\log N)^{1/2}}), & \text{unconditionally}, \\ (1-6\pi^{-1}e^{\gamma})\log\log N + O((\log N)^{-1/2}\log\log N), & \text{conditional on the RH}, \\ (1-6\pi^{-1}e^{\gamma})\log\log N + \Omega_{\pm}((\log N)^{-1/2}\log\log\log N/\log\log N), & \text{unconditional oscillations}, \end{cases}$

the nonsquarefree kernel of $\sigma(N)/N$.

Proof: Utilize the decomposition of the sum of divisors function

$$\sum_{d|N}\frac{1}{d} = \sum_{d|N,\,\mu(d)\neq 0}\frac{1}{d} + \sum_{d|N,\,\mu(d)=0}\frac{1}{d}, \tag{15}$$

and apply Theorem 4 to the identities

$$\sum_{d|N,\,\mu(d)\neq 0}\frac{1}{d} = \prod_{p\leq x}\left(1+\frac{1}{p}\right), \quad\text{and}\quad \sum_{d|N,\,\mu(d)=0}\frac{1}{d} \leq \left(\frac{\pi^2}{6}-1\right)\prod_{p\leq x}\left(1+\frac{1}{p}\right), \tag{16}$$

where $x = c_1 \log N$, and $c_1 > 0$ is a constant. ∎

The oscillating nature of the sum of divisor function is readily visible in the Ramanujan expansion

$$\frac{\sigma(N)}{N} = \frac{\pi^2}{6}\left(1 + \frac{\cos\pi N}{2^2} + \frac{2\cos 2\pi N/3}{3^2} + \cdots\right). \tag{17}$$

**2.2 The Extreme Values of Dedekind Psi Function.** The proof of Theorem 1 presented below relies on the oscillation theorem of the finite prime product; a completely elementary proof, not based on the oscillation theorem is also available, but it is longer.

From the elementary inequality

$$\prod_{p|N}(1+1/p) = \prod_{p|N}(1-1/p^2)\prod_{p|N}(1-1/p)^{-1} > \frac{6}{\pi^2}\prod_{p|N}(1-1/p)^{-1}, \tag{18}$$

and Theorems 3 and 4, it is plausible to expect that $\prod_{p|N}(1+1/p) > 6e^{\gamma}\pi^{-2}\log\log N$ infinitely often as $N \to \infty$.

***Theorem* 1.** Let $N \in \mathbb{N}$ be a primorial integer, then $\psi(N)/N > 6\pi^{-2}e^{\gamma}\log\log N$ holds unconditionally for all sufficiently large $N = 2\cdot 3\cdot 5\cdots p_k$.

Proof : Theorem 4 implies that the product

$$\prod_{p \leq x}(1+1/p) = \frac{6e^{\gamma}}{\pi^2}\log x + \Omega_{\pm}\left(\frac{\log\log\log x}{x^{1/2}\log\log x}\right). \tag{19}$$

In particular, it follows that

$$\prod_{p \leq x}(1+1/p) > \frac{6e^{\gamma}}{\pi^2}\log x + c_0\frac{\log\log\log x}{x^{1/2}\log\log x} \tag{20}$$

and

$$\prod_{p \leq x}(1+1/p) < \frac{6e^{\gamma}}{\pi^2}\log x - c_0\frac{\log\log\log x}{x^{1/2}\log\log x}$$

occur infinitely often as $x \to \infty$, where $c_0$, $c_1$, and $c_2 > 0$ are constants. It shows that $\prod_{p \leq x}(1+1/p)$ oscillates infinitely often, symmetrically about the line $6e^{\gamma}\pi^{-2}\log x$ as $x \to \infty$.

To rewrite the variable $x \geq 1$ in terms of the integer $N$, recall that the Chebychev function $\vartheta(x) = \sum_{p \leq x}\log p$, and

$$\log N_k = \sum_{p \leq p_k}\log p = \vartheta(p_k), \quad \text{and} \quad \vartheta(p_k) = p_k + o(p_k) \leq c_1\log N_k. \tag{21}$$

So it readily follows that $p_k \leq x = c_1\log N_k$.

Moreover, since the maxima of the sum of divisor function $\sigma(N) \geq \psi(N)$ occur at the colossally abundant integers $N = 2^{v_1} \cdot 3^{v_2} \cdots p_k^{v_k}$, and $v_1 \geq v_2 \geq \cdots \geq v_k \geq 1$, see [AE], [BR], [LA], [RJ], it follows that the maxima of the Dedekind psi function $\psi(N)$ occur at the squarefree primorial integers $N_k = 2\cdot3\cdot5\cdots p_k$. Therefore, expression (19) implies that

$$\begin{aligned}\frac{\psi(N_k)}{N_k} &= \prod_{p \leq c_0 \log N_k}(1+1/p) \\ &> \frac{6e^{\gamma}}{\pi^2}\log\log N_k + c_2\frac{\log\log\log N_k}{(\log N_k)^{1/2}\log\log N_k}\end{aligned} \tag{22}$$

as the primorial integer $N_k = 2\cdot3\cdot5\cdots p_k$ tends to infinity. ∎